\documentclass[3p,times]{elsarticle}

\usepackage{fullpage}
\usepackage{amssymb}
\usepackage{amsmath}
\usepackage{amsfonts}
\usepackage{subcaption}
\usepackage{multirow}
\usepackage{tikz}
\DeclareMathOperator\sign{sign}
\newtheorem{code}{Code}
\newtheorem{algorithm}{Algorithm}

\usepackage{xcolor}
\usepackage[colorlinks]{hyperref}
\definecolor{darkgreen}{rgb}{0.0,0.4,0.0}
\definecolor{darkred}{rgb}{0.6,0.0,0.0}
\definecolor{darkblue}{rgb}{0.0,0.0,0.5}
\definecolor{gray}{rgb}{0.5,0.5,0.5}
\definecolor{cyan}{rgb}{0.0,1.0,1.0}
\definecolor{darkcyan}{rgb}{0.0,0.5,0.5}
\definecolor{darkorange}{rgb}{0.8,0.4,0.0}
\definecolor{darkmargenta}{rgb}{0.5,0.0,0.5}
\definecolor{black}{rgb}{0.0,0.0,0.0}

\hypersetup{citecolor=blue, linkcolor=red, anchorcolor=darkblue,
filecolor=darkblue, urlcolor=darkblue}
\hypersetup{pdfauthor={Nicolae Suciu},pdftitle={Suciu}}

\journal{Advances in Water Resources}

\begin{document}

\begin{frontmatter}

\title{Programming with Chebfun. Case study: Richards equation}

\author{Nicolae Suciu\corref{1}} \ead{nsuciu@ictp.acad.ro}
\cortext[1]{Corresponding author.}
\address{Tiberiu Popoviciu Institute of Numerical Analysis, Romanian Academy, F\^{a}ntanele 57, 400320 Cluj-Napoca, Romania}

\begin{abstract}
The Chebfun software system is a Matlab extension essentially based on representations of (piece-wise) smooth one-variable functions by expansions in Chebyshev polynomials. One of Chebfun's attractive features is the ability to provide solutions to nonlinear boundary value problems (BVP) with accuracy close to the machine precision. This is done by the \texttt{chebop} class which provides automatic solutions by performing linearizations with a Newton method in function spaces of the nonlinear BVP, automatic differentiation, and using Fast Fourier Transform computations for the coefficients of the Chebyshev polynomials. A drawback of \texttt{chebop} automatic approach is the possible lack of convergence of the Newton method if the initial guess is not close enough to the exact solution. An explicit functional linearization done for each particular shape of the differential operator (i.e. without automatic differentiation) proves to be more robust than the \texttt{chebop} class and allows an enlargement of the range of convergence. Another alternative is the implicit $L$-scheme (quasi-Newton approach with derivatives replaced by suitable positive constants $L$), with a much simpler implementation and globally convergent. While \texttt{chebop} is the easiest way to solve the BVP, provided that it converges, the last two approaches largely overcome the convergence issue, yielding accurate solutions to a wide class of steady-state one-dimensional problems governed by Richards equation. Chebfun2 and Chebfun3, which at the current stage cannot solve BVPs, provide efficient tools for accuracy and convergence assessments of the unsteady solutions in one or two spatial dimensions obtained by classical discretization schemes.

\end{abstract}


\end{frontmatter}

\section{Introduction}
\label{sec:intro}

The Chebfun system is a free downloadable Matlab package enabling computing with functions similar to computing
with vectors aiming to ``feel symbolic but run at the speed of numerics''. While ``Chebfun'' with capital C denotes the software, ``chebfun'' is a function of one veriable defined on an interval of the real axis. The operations with chebfuns follow (almost) the same syntax used in operations with vectors in MATLAB, with commands for vectors overloaded by those acting on functions. The implementation of Chebfun is based on polynomial interpolation of smooth functions in Chebyshev points which, thanks to Fast Fourier Transform (FFT), finally results in expansions in Chebyshev polynomials \cite{Trefethen2019}.  

The basic theorem on Chebyshev series says that Lipschitz continuous functions have unique representations as Chebyshev series. That is, the series
$$f(x)=\sum_{k=0}^{\infty}a_kT_k(x),$$ where $T_k$ is the $k$-th order Chebyshev polynomial, is absolutely and uniformly convergent, and the coefficients $a_k$ have a closed form in terms of $T_k$ \cite[Theorem 3.1]{Trefethen2019}. 

Chebfun constructs approximations of a smooth function $f$ as follows. First, it calculates the polynomial interpolant 
$$f_n(x)\approx \sum_{j=0}^{n}a_jT_j(x),$$ at $n=9$ Chebyshev points $x_j=\cos(j\pi/n)$, $0\leq j \leq n$, and checks whether the Chebyshev coefficients are small enough. If this is not the case, the number of points is progressively increased until the norm $\|f-f_n\|$ falls to the level of rounding error. Further, the length of the series is reduced by neglecting the terms deemed to be negligible. Finally, the resulting series is evaluated at $n=$\texttt{length(f)} points via FFT to compute the Chebyshev coefficients defining the approximating chebfun \cite[Chap. 3 and 4]{Trefethen2019}.

\begin{figure}
\begin{minipage}[h]{0.45\linewidth}\centering
\includegraphics[width=\linewidth]{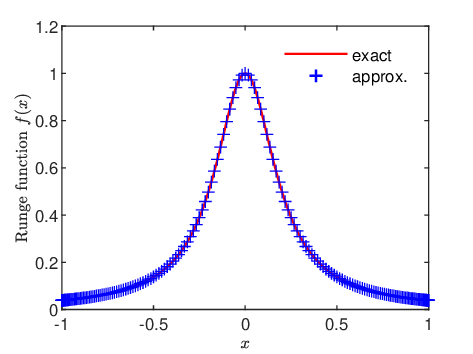}
\caption{\label{fig:runge}Chebfun approximation of the Runge function.}
\end{minipage}
\hspace*{0.1in}
\begin{minipage}[h]{0.45\linewidth}\centering
\includegraphics[width=\linewidth]{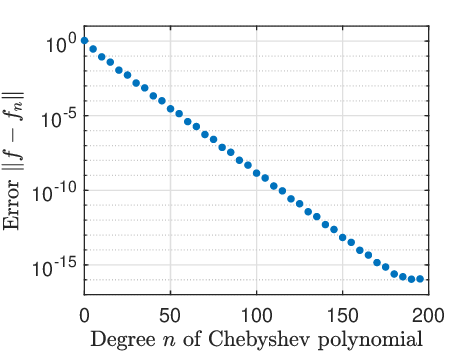}
\caption{\label{fig:geometric}Geometric decay of the error.}
\end{minipage}
\end{figure}

The accuracy of the approximation depends on the degree of smoothness of the function $f$. For instance, it $f$ is $\nu$-times differentiable, the approximation can achieve at most the algebraic accuracy $\|f-f_n\|=\mathcal{O}(n^{-\nu})$ \cite[Theorem 7.2]{Trefethen2019}. Instead, if (and only if) $f$ is (piece-wise) analytic one achieves the geometric accuracy 
$\mathcal{O}(C^{-n}), \;\; C>1$ \cite[Chap. 8]{Trefethen2019}. In the cited reference, this is illustrated by the chebfun approximation of the Runge function $f(x)=(1+25x^2)^{-1}$ by a Chebyshev polynomial of degree $n$ with the command  \texttt{fn = chebfun(f,n+1)}. Computing the approximation for $n=200$ shown in Fig.~\ref{fig:runge} takes about 0.5 seconds. The straight-line semi-log plot of the error shown in Fig.~\ref{fig:geometric} demonstrates the geometric accuracy of the chebfun approximation.

The geometric accuracy is not a special property of chebfun approximation, but also characterizes classical discretization methods. Finite element/volume/difference methods are usually characterized by algebraic error estimates, $\mathcal{O}\left(h^{\alpha \;} \right),\;h<1,\;\alpha >0$, with respect to the mesh size $h$. However, this is not similar to the algebraic accuracy with respect to $n$ of chebfun approximations for differentiable functions described above. For every mesh refinement $\{h_n=h_0r^{-n} \;| \;r>1\}$ the errors of the discretization methods can be recast as $\mathcal{O}(C^{- n}), \; C=r^{\alpha}>1$, which corresponds to the geometric accuracy.

The accuracy $\mathcal{O}(C^{- n})$ can be thought of as a functional dependence of the error $e_n=\|f-f_n\|$. It enables the assessment of the convergence order of the sequence $\{e_n\}$ according to the definition of the classical $C$-order $p\geq 1$ \cite{Catinas} by verifying the existence of the finite limit
\begin{equation}\label{eq:Cord}
Q_p=\lim\limits_{n\rightarrow\infty} \frac{e_{n+1}}{(e_n)^p}=\lim\limits_{n\rightarrow\infty} \frac{C^{-(n+1)}}{(C^{- n})^p}=\lim\limits_{n\rightarrow\infty} C^{n(p-1)-1}.
\end{equation}
We see that a finite limit $Q_p\in(0,\infty)$ can only be obtained if $p=1$. Since $C>1$, $Q_p=C^{-1}<1$, which corresponds to the linear convergence (as already indicated by the semi-log plot in Fig.~\ref{fig:geometric}). We have thus come to the somewhat surprising conclusion that both Chebfun and discretization approaches exhibit geometric accuracy and linear $C$-order of convergence. This could be misleading, because the mesh refinement becomes prohibitive from computational point of view after a few steps $n$. Instead, Chebfun effectively achieves the machine precision in remarkably short computing times.

An attractive feature of Chebfun is its ability to solve one-dimensional nonlinear boundary value problems (BVP). This is achieved with the object class \texttt{chebop} introduced since Chebfun version 4, which encompasses the class of linear operators as a particular case. Chebop uses forward automatic differentiation to compute Frech\'et derivatives, constructs linearized BVPs, and provides fully automatic solutions for a wide range of nonlinear BVPs \cite{BirkissonandDriscoll2012}. More specifically, \texttt{chebop} uses a variant of Newton's method implemented for continuous functions, rather than vectors, where the Jacobian matrices are replaced by Frech\'et derivative operators. Following \cite{BirkissonandDriscoll2012}, the procedure can be described as follows. Let's start with the the nonlinear BVP
\begin{align}
& \Phi(u)=0\label{eq:bvp1}\\
& \alpha(u)\vert_{x=a} = 0, \;\; \beta(u)\vert_{x=b} =0\label{eq:bvp2}
\end{align}
associated to the nonlinear differential operator $\Phi$ which solves for $u(x)$ in the interval $x\in[a, b]$ under the boundary conditions formulated with the aid of the (possibly nonlinear) differential operators $\alpha$ and $\beta$.
Generalizing the Taylor expansion used to derive the classical Newton method in finite dimensional spaces one starts with the current iterate $u_k$ and one evaluates $\Phi(u_k+v)$ for a given perturbation $v$ by the linearization $\Phi(u_k)+\Phi'(u_k)v$, where the prime denotes the Frech\'et derivative. The Frech\'et differential $\Phi'(u_k)v=\Phi(u_k+v)-\Phi(u_k)$ is computed with a standard perturbation approach by retaining linear terms in $v$ and dropping higher order terms. After replacing $\Phi(u)$ by $\Phi(u_k+v)$ in the differential equation (\ref{eq:bvp1}) and proceeding in the same way with the boundary conditions (\ref{eq:bvp2}) one obtains the update (or correction) $v$ as the solution of the linearized problem
\begin{align}
& \Phi'(u_k)v+\Phi(u_k)=0\label{eq:lin1}\\
& \alpha'(u_k)v\vert_{x=a}+\alpha(u_k)\vert_{x=a} = 0, \;\; \beta'(u_k)v\vert_{x=b}+\beta(u_k)\vert_{x=b} =0\label{eq:lin2}.
\end{align}
By solving (\ref{eq:lin1}-\ref{eq:lin2}) one obtains the new iterative solution to (\ref{eq:bvp1}-\ref{eq:bvp2}) given by $u_{k+1}=u_k+v$. The Newton iteration stops if any of the following criteria is satisfied,
\begin{align}
&\|v\|<\varepsilon \label{eq:crit1}\\
&\frac{1}{|u_k\|}\sqrt{\|\Phi(u_k)\|^2+\|\alpha(u_k)\vert_{x=a}\|^2+\|\beta(u_k)\vert_{x=a}\|^2}<\varepsilon \label{eq:crit2}\\
& \mbox{\texttt{iter} > \texttt{maxIter}}, \label{eq:crit3}
\end{align}
where $\varepsilon$ is a given tolerance, \texttt{iter} is the current iteration count, and \texttt{maxIter} is the user defined maximum number of iterations. 

The class \texttt{chebop} solves the linearized problem (\ref{eq:lin1}-\ref{eq:crit3}) by using automatic differentiation. In \cite{BirkissonandDriscoll2012}, the authors illustrate the approach by solving the nonlinear BVP (\ref{eq:bvp1}-\ref{eq:bvp2}) with $\Phi(u)=u''+2u\sin u, \; \alpha(u)=u', \; \beta(u)=uu'-2$. Carrying out the computations with the Chebfun version 4, in use when they published the paper, the authors reported identical results for both the explicit linearization (\ref{eq:lin1}-\ref{eq:lin2}) and the automated \texttt{chebop} approach implemented as in the sequence of commands presented below.
\begin{code}\label{code:chebop}
\texttt{chebop} class
\begin{verbatim}
N = chebop(0,5);                   % Create a chebop on [0,5]
N.op = @(u) diff(u,2)+2*u.*sin(u); % Assign the DE part
N.lbc = 'neumann';                 % Homogeneous Neumann LBC
N.rbc = @(u) diff(u).*u-2;         % Set RBC, u’(5)*u(5) = 2
N.init = chebfun('x',[0 5]);       % Initial guess of the solution
u = N\0;                           % Solve using overloaded \
\end{verbatim}
\end{code}
Repeating the computations with the current Chebfun version 5.7, we find that the explicit linearization works while the automated \texttt{chebop} fails to converge. We also find the if the definition interval is reduced from $[0,5]$ to $[0, 0.5]$ Newton's method converges in 6 iterations. When choosing in Chebfun preferences \cite[Chap. 8]{Drisscolletal2014} the option \texttt{cheboppref.setDefaults('display','iter')}, the output contains details about the iteration process and the ``Final error estimate: 4.53e-15 (differential equation)''. At this stage, we do not know whether the final error corresponds to the update (\ref{eq:crit1}) or to the residual (\ref{eq:crit2}). For more details, we have to solve with \texttt{[u,info] = solvebvp(N,0)} which provides the field \texttt{info.normDelta} containing the norms of the update $v$ at each iteration. In this case it is found that the final error estimate is precisely the norm (\ref{eq:crit1}) of the last iteration update. 

Comparing the two solution approaches, we can see that the explicit linearization is more robust than the automated \texttt{chebop} (actually, we found that it works for all the BVP problems considered in this study). Moreover, the explicit approach is more flexible. For instance, we can choose the stopping criterion, e.g. in terms of residuals of the differential operator. We also can easily specify the output: residuals of the solution $u$ and of the update $v$, errors with respect to a known solution, coefficients of the Chebyshev expansion, and so on. However, the method is constrained by the choice of a suitable initial solution. For the BVP discussed here, replacing the initial guess $u=x$ by $u=0.9x$ leads to the failure of the convergence.

To cope with the convergence issues, we propose a new linearization approach called $L$-scheme. This consists of a formulation for continuous functions of the implicit $L$-scheme used in classical discretization approaches. The latter is a stabilized Picard iteration \cite{Celiaetal1990} where the Jacobian is replaced by a suitably chosen positive parameter $L$  \cite{ListandRadu2016}. The $L$-scheme has two nice features. Fist, it does not require prior manual calculation of the terms of the linearized problem (\ref{eq:lin1}-\ref{eq:lin2}). It is enough to add the stabilization term $Lv$ and to perform an explicit linearization by expressing the nonlinear terms of the BVP operator with the solution of the previous iteration. The second feature is the global convergence of the scheme, as it only depends on the choice of the stabilization parameter $L$.

In this paper, we illustrate the use of the automated \texttt{chebop}, explicit linearization, and $L$-scheme by solving one-dimensional steady-state problems for infiltration in unsaturated soils governed by Richards equation. Further, we explore the capability of the newly coming Chebfun2 and Chebfun3 systems to solve unsteady infiltration problems 
in one and two spatial dimensions. This is done with a two-way coupling of Chebfun and finite difference (FD) codes: Chebfun codes provide exact manufactured solutions and source terms to FD $L$-schemes functions, the latter compute numerical solutions at fixed time and iteration steps, which are then imported in Chebfun codes, converted into chebfun, chebfun2, and chebfun3 objects, and further used to estimate errors and residuals. 
The last application of programming with Chebfun presented in this paper uses a similar coupling to evaluate orders of accuracy of the FD $L$-schemes.

The paper is organized as follows. Section~\ref{sec:lin} presents the principles and compares several linearization approaches that will be further used in section~\ref{sec:infil} to implement Chebfun solutions to BVP for steady state infiltration. Section~\ref{sec:cheb23} contains applications of Chebfun2 and Chebfun3 for coupled Chebfun-FD schemes and orders of accuracy estimates for FD $L$-schemes. Finally, some conclusions are drawn in section~\ref{sec:concl}. The Matlab codes, with Chebfun extension, used in this study are posted online at \href{https://github.com/PMFlow/programming\_with\_chebfun}{https://github.com/PMFlow/programming\_with\_chebfun}. References in the text are done by the folder and code names in brackets, [Folder/code.m].

\section{Linearization approaches}\label{sec:lin}

Newton's method is a prototypical linearization approach for nonlinear problems. Since the $L$-schemes can be thought of as quasi-Newton methods, we shall introduce them in the following by simplifying the Newton method while using its fixed-point structure. 

Let us consider the problem of finding the zeros of the one variable function $f(x)$. Looking for an iterative method to solve $f(x)=0$, we denote by $x_k$ the current approximation and solve the linearization obtained by the first-order Taylor expansion $f(x_{k+1})=f(x_k)+f'(x_k)(x_{k+1}-x_k))=0$ to obtain a new approximate solution,
\begin{equation}\label{eq:Nscheme}
x_{k+1}=x_k-\frac{f(x_k)}{f'(x_k)}=F(x_k).
\end{equation}

Hence, the scheme is equivalent to the fixed point iteration $x_{k+1}=F(x_k)$ and the solution $x^*$ is a fixed-point of $F$, $x^*=f(x^*)$. If $F$ is $p$-times differentiable the convergence of the fixed-point iterations is ensured under the following conditions \cite[Corollary 8.8]{KnabnerandAngermann2003}: the iteration $x_{k+1}=F(x_k)$ is locally convergent in a neighborhood $U(x^*)$ of the fixed point

$\bullet$ $if\;\; p=1,\;\; F'(x^*)\neq 0,\;\; and\;\; \vert F'(x^*)\vert<1\;$, which implies
\begin{align}\label{eq:p1}
\vert x_{k+1}-x^*\vert\leq C\vert x_k-x^*\vert,\;\; C<1,
\end{align}  

$\bullet$ $if\;\; p>1,\;\; F'(x^*)=F''(x^*)\ldots F^{(p-1)}(x^*)=0, \;\; and\;\; F^{(p)}\neq 0,\;$ which implies
\begin{align}\label{eq:pgt1}
\vert x_{k+1}-x^*\vert\leq C\vert x_k-x^*\vert^p,\;\; C>0.
\end{align}

An explicit $L$-scheme may be constructed as a relaxation of the Newton scheme, with $f'(x_k)$ in (\ref{eq:Nscheme}) replaced by a positive constant $L$. Let $f\in C^1(\mathbb{R})$ with a unique root $x^*$. Choosing $L>\sup_{x\in\mathbb{R}}f'(x)$, we define the explicit $L$-scheme by
\begin{equation}\label{eq:Lscheme}
x_{k+1}=x_k-\sign(f'(x_k))\frac{f(x_k)}{L}=F_L(x_k).
\end{equation}
The sign function in (\ref{eq:Lscheme}) ensures the condition 
$$\vert F_{L}'(x_k)\vert=\left\vert 1-\sign(f'(x_k))\frac{f'(x_k)}{L}\right\vert<1$$ for all $x\in\mathbb{R}$ 
regardless of whether the derivative $f'(x_k)$ is positive or negative. Thus, according to (\ref{eq:p1}) the $L$-scheme converges linearly to $x^*$. Since the condition $\vert F_{L}'(x^*)\vert<1$ is not restricted to a neighborhood of $x^*$, like for (\ref{eq:p1}), the $L$-scheme is globally convergent.

As an illustration let us consider the equation \cite[Example 4.6]{Catinas}
\begin{align}\label{eq:fx}
f(x)=x+x^2+10x^3=0,\;\; with\;fixed\;point\;\; x^*=0
\end{align}  
and the corresponding fixed-point function $F(x)=x-f(x)/f'(x)$. Since $F'(x^*)=0$ and $F''(x^*)=2$ the fixed point iteration converges and behaves as (\ref{eq:pgt1}) with $p=2$. Moreover, it is shown \cite[Theorem 4.2]{Catinas} that the sequence of errors converges with the classical $C$-order $p=2$, according to (\ref{eq:Cord}). 

\begin{figure}
\begin{minipage}[h]{0.45\linewidth}\centering
\includegraphics[width=\linewidth]{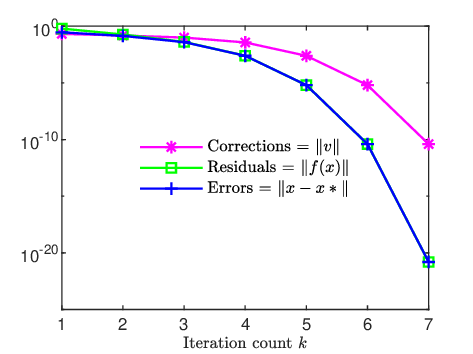}
\caption{\label{fig:Newton}Convergence of the Newton scheme.}
\end{minipage}
\hspace*{0.1in}
\begin{minipage}[h]{0.45\linewidth}\centering
\includegraphics[width=\linewidth]{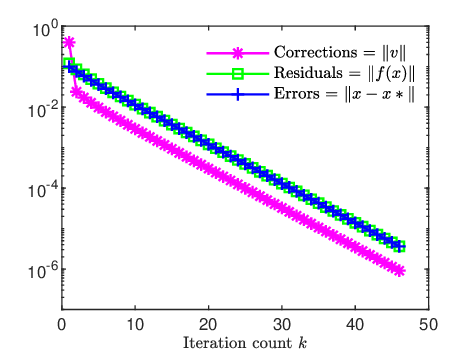}
\caption{\label{fig:Lscheme}Convergence of the $L$-scheme.}
\end{minipage}
\end{figure}

The convergence of the Newton scheme (\ref{eq:Nscheme}) solving the problem $(\ref{eq:fx})$ with the code [\href{https://github.com/PMFlow/programming_with_chebfun/blob/main/Linearization/Newton_Fx.m}{Linearization/Newton\_Fx}] is quantified by $L_2$ norms of corrections, $v=(x_{k+1}-x_k)$, errors, $(x_{k+1}-x^*)$, and residuals, $f(x_{k+1})$. With the initial guess $x=0.5$ and stopping criterion for corrections given by $\|v\|<10^{-6}$, the convergence is achieved in 7 iterations. The concave shape of the plots shown in Fig.~\ref{fig:Newton} indicates a convergence order $p>1$ \cite{Catinas}. The $L$-scheme, with parameter $L=5$, converges in 46 iterations and Fig.~\ref{fig:Lscheme} indicates the linear convergence. Note that when no criteria for choosing the parameter $L$ are available, we progressively increase it until the convergence occurs.

In a space of continuously differentiable functions $u(x)$, the analogue of the fixed-point problem (\ref{eq:fx}) reads
\begin{align}\label{eq:fx_cheb}
\Phi(u)=u+u^2+10u^3=0,\;\; with\;\; fixed\;\; point \;\; u^*\equiv0.
\end{align}  
The problem (\ref{eq:fx_cheb}) can be solved in Chebfun by linearization procedures similar to those presented above with the commands [\href{https://github.com/PMFlow/programming_with_chebfun/blob/main/Linearization/Frechet_Fx.m}{Linearization/Frechet\_Fx}]:

\begin{code}\label{code:iter}
Linearization methods
\begin{verbatim}
D=[-1 1];                                         % domain
x = chebfun('x',D);                               % chebfun 'x' on D
A = chebop(D);                                    % linear operator on D
u_ref = chebfun(0,D);                             % exact solution 
u0=chebfun('5*sin(pi*x/2',D);                     % initial guess
phi = @(u) u+u.^2+10*u.^3;                        % anonymous function for Eq. (14)
dF = @(u) 30*u.^2+2*u+1;                          % derivative of 'phi'
u=u0;                                             % initial iterate
nrmv = 1;               
while nrmv > 1e-6
    if iop==1
        L=200
        u=u0-sign(dF(u))*phi(u0)/L;               % explicit L-scheme 
        v = u-u0;
    elseif iop==2
        L=5;
        A.op = @(u) L*(u-u0)+sign(dF(u0)).*(u+u0.*u+10*u0.^2.*u); % implicit L-scheme 
        u=A\0;  
        v = u-u0;
    elseif iop==3
        A.op = @(v) (30*u.^2+2*u+1)*v+phi(u);     % Newton-Frechet in function space
        v = A\0;                                   
        u = u+v;
    end
    nrmv = norm(v); 
    u0=u;
end
\end{verbatim}
\end{code}

The explicit $L$-scheme in (\ref{code:iter}) is an extension to function space of the scheme (\ref{eq:Lscheme})  for finite-dimensional spaces. The implicit $L$-scheme, inspired by approaches from discretization methods \cite{ListandRadu2016} consists of adding the stabilization term $L(u-u_0)$ to the equation (\ref{eq:fx_cheb}) linearized by using the solution $u_0$ of the previous iteration, $u+u_0u+10u_{0}^{2}u$. Unlike the scheme presented in \cite{ListandRadu2016}, which exploits the monotonic increase of the nonlinear term (derivative of the saturation) we use the factor $\sign(\Phi')$, as in the explicit scheme (\ref{eq:Lscheme}), to cope with possible sign changes of the slope of $\Phi$. The third scheme in the Code~(\ref{code:iter}) corresponds to the Newton-Frec\'echet linearization (\ref{eq:lin1}) underlying the automatic \texttt{chebop} approach.

\begin{figure}
\begin{minipage}[h]{0.45\linewidth}\centering
\includegraphics[width=\linewidth]{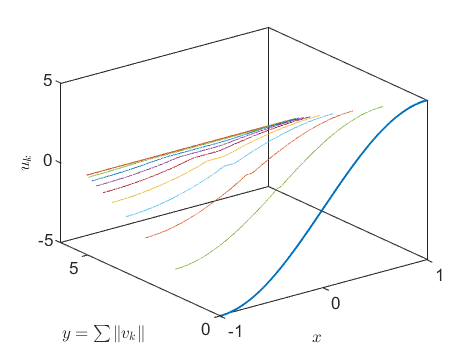}
\caption{\label{fig:Frechet_sol}Successive iterative solutions of the Fr\'echet scheme.}
\end{minipage}
\hspace*{0.1in}
\begin{minipage}[h]{0.45\linewidth}\centering
\includegraphics[width=\linewidth]{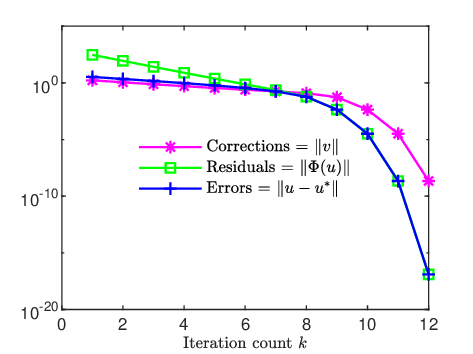}
\caption{\label{fig:Frechet_conv}Convergence of the Fr\'echet scheme.}
\end{minipage}
\end{figure}

With an initial guess given by the constant function \texttt{u0=chebfun(0.5,D)}, similar to the finite dimensional case above, the results of Frec\'echet linearization scheme are (almost) identical to those of the Newton scheme presented in 
Fig.~\ref{fig:Newton} and the same holds for the results of the two $L$-schemes in function space and those from Fig.~\ref{fig:Lscheme}. These results are what we expected, because the constant fixed-point solution $u=0$ is approached by constant functions, starting with $u_0=0.5$, and provide a validation of the linearisation schemes in function space used in Code~\ref{code:iter}. The three schemes exhibit different performance when the initial guess is non-trivial. With \texttt{u0=chebfun('5*sin(pi*x/2',D)}, the explicit $L$-scheme takes about 90 seconds to achieve the convergence, after 1500 iterations, with final errors and residuals of order $10^{-4}$. (All the computations presented in this paper are done on a laptop with 16 GB installed RAM and 1.9 GHz.) The implicit $L$-scheme convergence after about 30 seconds and 60 iterations, with final errors and residuals of order $10^{-6}$. The convergence is ensured by taking $L=5$ in the implicit scheme, while the explicit scheme requires a much higher parameter $L=200$. The Fr\'echet linearization proves to be the most performant in solving this problem. It converges in about 9 seconds, after 12 iterations, with final errors and residuals of order $10^{-17}$. The approach of the space-variable initial iterate to the constant fixed-point solution is shown in Fig.~\ref{fig:Frechet_sol} and the behavior of the corrections, residuals, and errors with respect to the iteration count is presented in Fig~\ref{fig:Frechet_conv}. 

\section{Application to steady-state infiltration in soils}\label{sec:infil}

Steady-state infiltration problems modeled by Richards equation provide a methodological framework for the determination of soil hydraulic properties and in advancing the knowledge of the soil physics \cite{LogsdonandJaynes1993,SiandKachanoski2000,BassetandAbouNajm2025}. Such investigations use numerical \cite{Suciuetal2021,Xiaoetal2022,Autovinoetal2024} as well as approximate analytical solutions to solve the nonlinear Richards equation \cite{Basha1999,BraddockandParlange2000,BraddockandParlange2003,TalukdarandBarua2022}. Since Chebfun system provides solutions to BVP problems with practically the machine precision, it could be a valuable alternative in solving steady-state infiltration problems.

\subsection{Problem formulation}
\label{sec:formulation}

Infiltration in soils is often modeled by the one-dimensional Richards equation
\begin{equation}\label{eq:Richards}
\frac{\partial}{\partial t}\theta(\psi)-\frac{\partial}{\partial z}\left[K(\psi)\frac{\partial}{\partial z}(\psi-z)\right]=f,
\end{equation}
where $\psi(z,t)$ is the pressure head expressed in length units, $\theta$ is the volumetric water content, $K$ stands for the hydraulic conductivity of the medium, $z$ is the height oriented positively downward, and $f$ is a sink/source term. The steady-state unsaturated flow occurs when the water content $\theta$ is still smaller than the porosity, possibly spatially variable, but time independent, $\tfrac{\partial}{\partial t}\theta(\psi)=0$. In the steady-state case of one-dimensional unsaturated flows it is convenient to work with the pressure potential defined by the pressure head with opposite sign \cite{Gardner1958}. With these, the steady-state Richards equation (\ref{eq:Richards}) becomes
\begin{equation}\label{eq:uptake}
\frac{d}{d z}\left[K(\psi)\left(\frac{d\psi}{d z}+1\right)\right]=-f,
\end{equation}
where $f$ is a positive function modelling the root water uptake by plants. Note that when working with the pressure potential $-\psi$ the water flux $q=-K(\frac{\partial \psi}{\partial z}-1)$, expressed in terms of pressure head $\psi$, becomes $q=K(\frac{\partial \psi}{\partial z}+1)$.

To close equation (\ref{eq:uptake}) we consider two model relationships between hydraulic conductivity and potential proposed by Gardner \cite{Gardner1958},
\begin{align}
&K(\psi)=a\exp(-c\psi),\label{eq:K1}\\
&K(\psi)=1/(1+\psi^n). \label{eq:K2}
\end{align}
where $a,\; c,\; \mbox{and}\; n$ are model constants. The closure (\ref{eq:K1}) is the classical exponential Gardner's closure further referred to as Gardner model. The second closure, (\ref{eq:K2}), has been used in a couple of papers by Basha, Braddock, and Parlange \cite{Basha1999,BraddockandParlange2000,BraddockandParlange2003} to derive exact and approximate solutions of the equation (\ref{eq:uptake}) and will be further referred to as Basha model. 

\subsection{Implicit $L$-scheme for Richards Equation}\label{sec:L-Richards}

The $L$-scheme for Richards equation can be obtained as a relaxation of the modified Picard method proposed by \cite{Celiaetal1990}. The latter is derived from the backward Euler time discretization of (\ref{eq:Richards}) and a simple Picard iteration scheme, with iteration count denoted by $k$,
\[
\frac{\theta^{j,k+1}-\theta^{j-1}}{\tau}-\frac{\partial}{\partial z}\left[K^{j,k}\frac{\partial}{\partial z}(\psi^{j,k+1}-z)\right]=0,
\]
followed by a truncated Taylor expansion of $\theta^{j,k+1}$ with respect to $\psi$ about $\psi^{j,k}$,
$$\theta^{j,k+1}=\theta^{j,k}+\left.\frac{d\theta}{d\psi}\right\vert^{j-1,k}(\psi^{j,k+1}-\psi^{j,k}),$$
and ends up to a modified Picard scheme
\begin{equation}\label{eq:modPicard}
\left.\frac{d\theta}{d\psi}\right\vert^{j-1,k}(\psi^{j,k+1}-\psi^{j,k})+\frac{\theta^{j,k}-\theta^{j-1}}{\tau}
-\frac{\partial}{\partial z}\left[K^{j,k}\frac{\partial}{\partial z}(\psi^{j,k+1}-z)\right]=0.
\end{equation}
The relaxation consists of replacing the first term in (\ref{eq:modPicard}) by the stabilization term 
$L(\psi^{j,k+1}-\psi^{j,k}),\; L>0$. In the steady-state case, the implicit $L$-scheme iteration reads
\begin{equation}\label{eq:L-Richards}
L(\psi^{j,k+1}-\psi^{j,k})-\frac{\partial}{\partial z}\left[K^{j,k}\frac{\partial}{\partial z}(\psi^{j,k+1}-z)\right]=0.
\end{equation}

The convergence of the $L$-scheme (\ref{eq:L-Richards}) has been proven in a variational formulation by exploiting the monotonicity of $\theta$ and assuming common regularity conditions on $K$ and the gradient of the solution \cite{ListandRadu2016}.

\subsection{Chebfun solutions with Gardner model}\label{sec:Gardner}

In the following we look for solutions of the equation (\ref{eq:uptake}) for $f=0$, in the domain $D=[0\; 2]$. A constant flux $q=0.1$ is imposed on the left boundary and the Dirichlet condition $\psi=1$ on the right boundary. The problem is closed by the Gardner model (\ref{eq:K1}) with $a=1$ and $c\in\{0.5,\; 1,\; 1.5,\; 2\}$. 

After inserting (\ref{eq:K1}) into (\ref{eq:uptake}), using the above boundary conditions, and denoting the potential $\psi$ by $u$ the nonlinear BVP to be solved reads
\begin{align}
& \Phi(u)=u''-cu'(u'+1)=0,\label{eq:bvpG1}\\
& \alpha(u)\vert_{x=0}=a[u'(0)+1]-q\exp[cu(0)]=0,\label{eq:bvpG2}\\
& \beta(u)\vert_{x=2}=u(2)-1=0.\label{eq:bvpG3}
\end{align}

The Chebfun solver based on \texttt{chebop} class, similar to code~\ref{code:chebop}, is given by the following code lines [\href{https://github.com/PMFlow/programming_with_chebfun/blob/main/Linearization/Chebop_Gardner_Fig1.m}{Linearization/Chebop\_Gardner\_Fig1.m}]:
\begin{verbatim}
N.op = @(x,u) u*diff(u,2)-c*diff(u)*(diff(u)+1);
N.lbc=@(u) a*(diff(u)+1)-q*exp(c*u);
N.rbc=@(u) u-b;
u = N\0;
\end{verbatim}
For $c=1.5$ and $c=2$ we got the message ``Warning: Newton iteration failed. Please try supplying a better initial guess via the .init field of the chebop.'' This is a common issue since it is known that the default \texttt{N.init}, constructed as the lowest-degree polynomial satisfying the boundary conditions \cite[Chap. 10]{Drisscolletal2014}, encounters difficulties to satisfy Newmann conditions \cite{BirkissonandDriscoll2012}. Since the issue could not be solved by trying different \texttt{N.init} fields we conclude that the results obtained with the \texttt{chebop} class are not reliable and proceed with other two alternative solution approaches.

\begin{figure}
\begin{minipage}[h]{0.45\linewidth}\centering
\includegraphics[width=\linewidth]{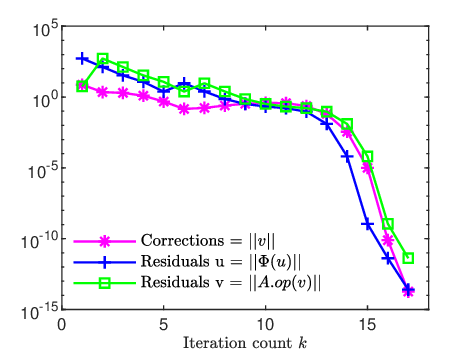}
\caption{\label{fig:errFrechetGardnerFig1}Convergence of the Fr\'echet scheme for BVP (\ref{eq:bvpG1}-\ref{eq:bvpG3}) with Gardner model, $f=0$, $a=1$ and $c=2$.}
\end{minipage}
\hspace*{0.1in}
\begin{minipage}[h]{0.45\linewidth}\centering
\includegraphics[width=\linewidth]{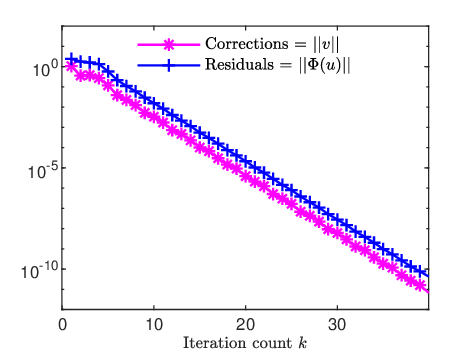}
\caption{\label{fig:errLschemeGardnerFig1}Convergence of the $L$-scheme for BVP (\ref{eq:bvpG1}-\ref{eq:bvpG3}) with Gardner model, $f=0$, $a=1$ and $c=2$.}
\end{minipage}
\end{figure}

In the following, we construct the explicit Newton-Fr\'echet linearization (\ref{eq:lin1}-\ref{eq:lin2}) of the problem (\ref{eq:bvpG1}-\ref{eq:bvpG3}). First, we compute the Fr\'echet differential of the nonlinear operator (\ref{eq:bvpG1}) for a perturbation $v$ about $u$,
\begin{align*}
d\Phi&=(u+v)''-c(u+v)'[(u+v)'+1]-u''+cu'(u'+1)\\
&=v''-c(2u'+1)v'.
\end{align*}
In the same way, we compute the Fr\'echet differential of (\ref{eq:bvpG2}),
\begin{align*}
d\alpha&=a(u'+v'+1)-q\exp[c(u+v)]-a(u'+1)+q\exp(cu)\\
&=av'-q\exp(cu)[\exp(cv)-1]\\
&=av'-q\exp(cu)cv,
\end{align*}
where we use the linearization $\exp(cv)=1+cv$.
Finally, for the right boundary condition (\ref{eq:bvpG3}) we have
$$d\beta=u+v-1-u+1=v.$$
According to (\ref{eq:lin1}-\ref{eq:lin2}), the Newton correction $v_k$ is obtained by solving the linearization problem
\begin{align*}
&v''_k-c(2u'+1)v'_k+u''-cu'(u'+1)=0,\\
&av'_k(0)-q\exp(cu(0))cv_k(0)+a[u'(0)+1]-q\exp[cu(0)]=0,\\
&v_k+u(2)-1=0.
\end{align*}
The Newton-Fr\'echet method is implemented in Chebfun [\href{https://github.com/PMFlow/programming_with_chebfun/blob/main/Linearization/Lscheme_Gardner_Fig1.m}{Linearization/Newton\_Frechet\_Gardner\_Fig1.m}] by
\begin{verbatim}
A.op = @(v) diff(v,2)-c*(2*diff(u)+1)*diff(v);                           
Du = diff(u);                                                            
A.lbc = @(v) a*diff(v)-q*exp(c*u(x1))*(c*v)+a*(Du(x1)+1)-q*exp(c*u(x1)); 
A.rbc = @(v) v+u(x2)-1;                                                  
v = A\0;
u = u+v;                                                         
\end{verbatim}
With initial guess $u=x$ and tolerance $\|v\|<10^{-11}$, the Fr\'echet scheme converges in about 17 seconds. The behavior of the corrections, and of the residuals for the solution $u$ and the update $v$ for the case $c=2$ are shown in Fig.~\ref{fig:errFrechetGardnerFig1}. 

\begin{figure}
\begin{minipage}[h]{0.45\linewidth}\centering
\includegraphics[width=\linewidth]{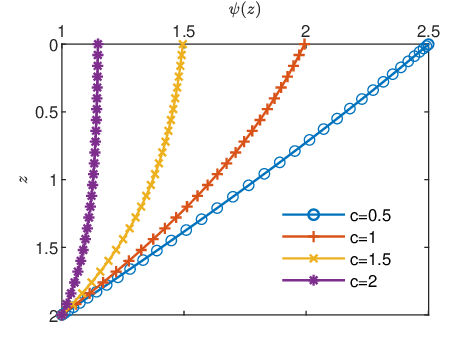}
\caption{\label{fig:GardnerFig1}Solution of the BVP (\ref{eq:bvpG1}-\ref{eq:bvpG3}) with Gardner model, $f=0$, and different parameters $c$.}
\end{minipage}
\hspace*{0.1in}
\begin{minipage}[h]{0.45\linewidth}\centering
\includegraphics[width=\linewidth]{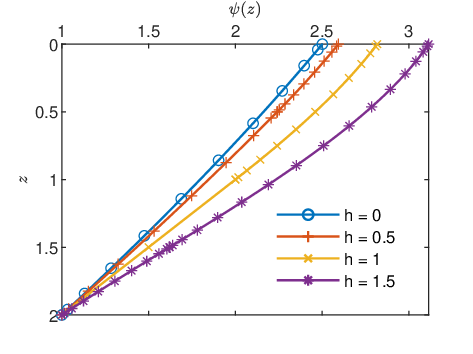}
\caption{\label{fig:GardnerFig2}Solution of the BVP (\ref{eq:bvpG1}-\ref{eq:bvpG3}) with Gardner model, $f>0$, and different uptake depths $h$.}
\end{minipage}
\end{figure}

The third solution approach used to solve (\ref{eq:bvpG1}-\ref{eq:bvpG3}) is the implicit $L$-scheme. By analogy with (\ref{eq:L-Richards}) we propose the scheme
\begin{equation}\label{eq:implicitL}
L(u-u_0)-Op(u,u_0)=0,
\end{equation}
where $u$ and $u_0$ are the current and previous iterates and $Op$ is a problem specific linear operator corresponding to (\ref{eq:uptake}). The $L$-scheme is implemented in Chebfun [\href{https://github.com/PMFlow/programming_with_chebfun/blob/main/Linearization/Lscheme_Gardner_Fig1.m}{Linearization/Lscheme\_Gardner\_Fig1.m}] by:
\begin{verbatim}
A.op  = @(u) L*(u-u0)-(diff(u,2)-c*diff(u0)*(diff(u)+1));  
A.lbc = @(u) a*(diff(u)+1)-q*exp(c*u0);                    
A.rbc = @(u) u-b;                                          
u = A\0;                                                   
\end{verbatim}
The $L$-scheme in function space, with the same initial guess and tolerance as above, converges in about 25 seconds. The linear convergence is proven by the decay of the corrections below the tolerance $\varepsilon$ shown in Fig.~\ref{fig:errLschemeGardnerFig1}. Since the final residuum is of the order $10^{-10}$, the solution of the $L$-scheme is practical indistinguishable from the solution of the Fr\'echet scheme shown in Fig.~\ref{fig:GardnerFig1}. 

In the case of the non-homogeneous equation (\ref{eq:uptake}) with $f>0$, describing the water uptake, (\ref{eq:bvpG1}) has to be replaced by $\Phi(u)=u''-cu'(u'+1)+\exp(cu)f/a=0$. The boundary conditions (\ref{eq:bvpG2}-\ref{eq:bvpG3}) do not change because they are independent of $f$. The water uptake is taken as a piece-wise constant function, $f(x)=qh$ for $x\leq h$, $f(x)=0$ for $x>h$. The problem is solved for $c=0.5$ and  $h\in\{0,\; 0.5,\; 1,\; 1.5\}$. 

The \texttt{chebop} code [\href{https://github.com/PMFlow/programming_with_chebfun/blob/main/Linearization/Chebop_Gardner_Fig2.m}{Linearization/Chebop\_Gardner\_Fig2.m}] does not converge with the default initial guess \texttt{N.init} but it converges for suitable constant initial guess. For \texttt{N.init=0} the convergence is achieved in about 7 seconds. The solutions for the four uptake depths $h$ are shown in Fig.~\ref{fig:GardnerFig2}. Their accuracy, evaluated by the residuals (\ref{eq:crit2}), ranges between $10^{-12}$ and $10^{-10}$.

The code [\href{https://github.com/PMFlow/programming_with_chebfun/blob/main/Linearization/Newton_Frechet_Gardner_Fig2.m}{Linearization/Newton\_Frechet\_Gardner\_Fig2.m}] converges in about 14 seconds, with residuals of the order $10^{-15}$ for all $h$ values, for constant \texttt{N.init} but it fails to converge for non-constant initial guess, e.g.~\texttt{N.init=x}. The $L$-scheme code [\href{https://github.com/PMFlow/programming_with_chebfun/blob/main/Linearization/Lscheme_Gardner_Fig2.m}{Linearization/Lscheme\_Gardner\_Fig2.m}] converges in about 22 seconds with residuals of the order $10^{-11}$, being instead globally convergent, e.g.~\texttt{N.init=x} or constant chebfuns.

\subsection{Chebfun solutions with  Basha model}\label{sec:Basha}

With the Basha model (\ref{eq:K2}) inserted into (\ref{eq:uptake}) and the same boundary conditions as in Section~\ref{sec:Gardner}, the nonlinear BVP, in the general case $f\geq 0$, reads
\begin{align}
& \Phi(u)=(1+u^n)u''-nu^{(n-1)}u'(u'+1)+(1+u^n)^2f=0,\label{eq:bvpB1}\\
& \alpha(u)\vert_{x=0}=u'(0)-(1+u^n(0)q+1=0,\label{eq:bvpB2}\\
& \beta(u)\vert_{x=2}=u(2)-1=0.\label{eq:bvpB3}
\end{align}

For the homogeneous case, $f=0$, Basha \cite{Basha1999} derived an exact solution for integer parameter $n$. This solution has been further used to validate approximate solutions for the general non-homogeneous problem \cite{Basha1999,BraddockandParlange2000}. Here, instead of using exact solutions for validation purposes, we shall construct Chebfun solution with the machine precision by the three linearization approaches described in Section~\ref{sec:Gardner} for both the homogeneous and the non-homogeneous cases.

\begin{figure}
\begin{minipage}[h]{0.45\linewidth}\centering
\includegraphics[width=\linewidth]{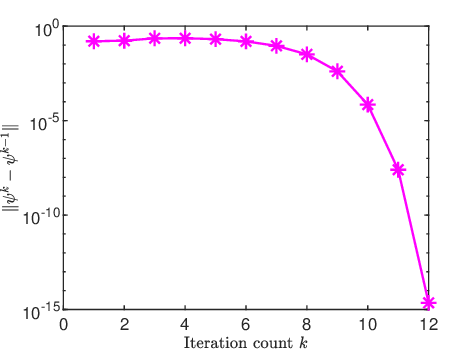}
\caption{\label{fig:conv_chebop}Convergence of the Chebop solution for BVP (\ref{eq:bvpB1}-\ref{eq:bvpB3}) with Basha model, $f=0$, $n=7$, and $q=0.1$.}
\end{minipage}
\hspace*{0.1in}
\begin{minipage}[h]{0.45\linewidth}\centering
\includegraphics[width=\linewidth]{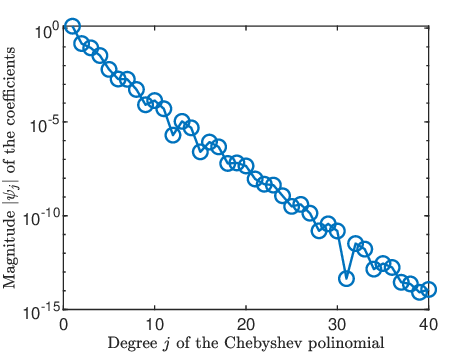}
\caption{\label{fig:coeff_chebop}Coefficients of the Chebyshev polynomial for BVP (\ref{eq:bvpG1}-\ref{eq:bvpG3}) with Basha model, $f=0$, $n=7$, and $q=0.1$.}
\end{minipage}
\end{figure}

\begin{figure}
\begin{minipage}[h]{0.45\linewidth}\centering
\includegraphics[width=\linewidth]{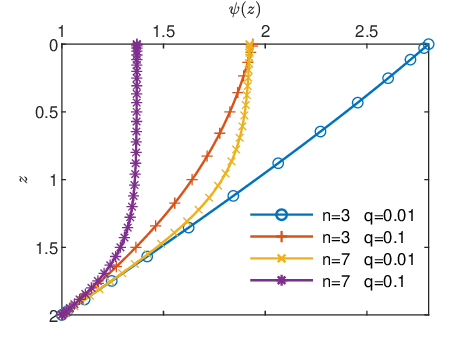}
\caption{\label{fig:BashaFig1}Solution of the BVP (\ref{eq:bvpG1}-\ref{eq:bvpG3}) with Basha model, $f=0$, and different parameters $n$ and $q$.}
\end{minipage}
\hspace*{0.1in}
\begin{minipage}[h]{0.45\linewidth}\centering
\includegraphics[width=\linewidth]{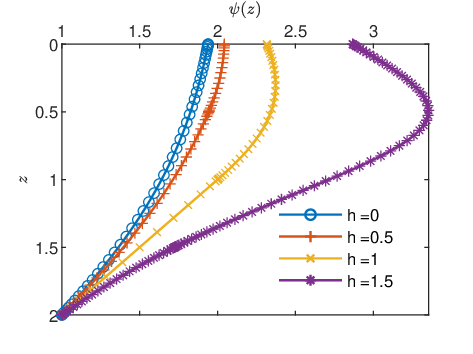}
\caption{\label{fig:BashaFig2}Solution of the BVP (\ref{eq:bvpG1}-\ref{eq:bvpG3}) with Basha model, $f>0$, and different uptake depths $h$.}
\end{minipage}
\end{figure}

For the beginning, we consider $f=0$ and the four combinations of parameters from \cite[Fig. 1]{Basha1999}: $\{n=3,\; q=0.01\}$, $\{n=3,\; q=0.1\}$, $\{n=7,\; q=0.1\}$, and $\{n=7,\; q=0.01\}$. The \texttt{chebop} code [\href{https://github.com/PMFlow/programming_with_chebfun/blob/main/Linearization/Chebop_Basha_Fig1.m}{Linearization/Chebop\_Basha\_Fig1.m}] uses a straightforward implementation with Chebfun syntax of the BVP (\ref{eq:bvpB1}-\ref{eq:bvpB3}). We are now in the happy situation where the convergence is achieved with the default initial guess in about 10 seconds. Figure~\ref{fig:conv_chebop} shows that the last correction coincides with the output ``Final error estimate: 2.25e-15 (differential equation)'', close to the machine precision. The coefficients of the Chebyshev polynomial (Fig.~\ref{fig:coeff_chebop}), which provide an indication about the approximation error, also approach the machine precision. Note that while the corrections decay following a concave curve consistent with a convergence order $p>1$, the coefficients follow a linear decay, corresponding to the convergence order $p=1$. Using (\ref{eq:crit2}) we also computed the final residual normalized to the last guess of the solution, $2.70\cdot 10^{-9}$, which provides a more reliable error estimate. The solution shown in Fig.~\ref{fig:BashaFig1} reproduces \cite[Fig. 1]{Basha1999}.

The Newton-Fr\'echet linearization method [\href{https://github.com/PMFlow/programming_with_chebfun/blob/main/Linearization/Newton_Frechet_Basha_Fig1.m}{Linearization/Newton\_Frechet\_Basha\_Fig1.m}] converges in about 12 seconds with residuals of the order $10^{-15}$. The $L$-scheme [\href{https://github.com/PMFlow/programming_with_chebfun/blob/main/Linearization/Lscheme_Basha_Fig1.m}{Linearization/Lscheme\_Basha\_Fig1.m}] converges in about 45 seconds with residuals between $10^{-8}$ and $10^{-9}$.

The non-homogeneous problem, with Basha model, is formulated with $n=3$ and $q=0.1$ for the same uptake depths as in Section~\ref{sec:Gardner}. The \texttt{chebop} implementation fails to converge with the default initial guess [\href{https://github.com/PMFlow/programming_with_chebfun/blob/main/Linearization/Chebop_Basha_Fig2.m}{Linearization//Chebop\_Basha\_Fig2.m}]. The convergence also fails for constant initial guess or initials guess given by a finite difference solution imported into the \texttt{chebop} code via different interpolation procedures.

The Newton-Fr\'echet code [\href{https://github.com/PMFlow/programming_with_chebfun/blob/main/Linearization/Newton_Frechet_Basha_Fig2.m}{Linearization/Newton\_Frechet\_Basha\_Fig2.m}] solves the cases shown in Fig.~\ref{fig:BashaFig2} in about 49 seconds. The convergence for $h>0$ is linear in this case, with residuals between $10^{-9}$ and $10^{-11}$. The $L$-scheme [\href{https://github.com/PMFlow/programming_with_chebfun/blob/main/Linearization/Lscheme_Basha_Fig2.m}{Linearization/Lscheme\_Basha\_Fig2.m}] converges linearly in about 80 seconds.

An overview of the results presented above is given in Table~\ref{tab:summary}.

\begin{table}[h]
\begin{tabular}{l|llll|llll}
\hline
\multirow{3}{*}{} & \multicolumn{4}{l|}{Gardner model}                                                    & \multicolumn{4}{l}{Basha model} \\ \cline{2-9} 
& \multicolumn{2}{l|}{$f=0$} & \multicolumn{2}{l|}{$f>0$} & \multicolumn{2}{l|}{$f=0$} & \multicolumn{2}{l}{$f>0$} \\ \cline{2-9} 
& \multicolumn{1}{l|}{conv.} & \multicolumn{1}{l|}{$\|\Psi(u)\|$} & \multicolumn{1}{l|}{conv.} &  \multicolumn{1}{l|}{$\|\Psi(u)\|$} & \multicolumn{1}{l|}{conv.} & \multicolumn{1}{l|}{$\|\Psi(u)\|$} & \multicolumn{1}{l|}{conv.} & \multicolumn{1}{l}{$\|\Psi(u)\|$} \\ \hline
\texttt{chebop} & \multicolumn{1}{l|}{NO} & \multicolumn{1}{l|}{--} & \multicolumn{1}{l|}{YES} &  \multicolumn{1}{l|}{6.33e-10} & \multicolumn{1}{l|}{YES} & \multicolumn{1}{l|}{5.00e-09} & \multicolumn{1}{l|}{NO} & \multicolumn{1}{l}{--} \\ \hline
Fr\'echet & \multicolumn{1}{l|}{YES} & \multicolumn{1}{l|}{2.52e-14} & \multicolumn{1}{l|}{YES} &  \multicolumn{1}{l|}{2.53e-16} & \multicolumn{1}{l|}{YES} & \multicolumn{1}{l|}{5.43e-15} & \multicolumn{1}{l|}{YES} & \multicolumn{1}{l}{3.67e-9} \\ \hline
$L$-scheme & \multicolumn{1}{l|}{YES} & \multicolumn{1}{l|}{4.16e-11} & \multicolumn{1}{l|}{YES} &  \multicolumn{1}{l|}{2.59e-11} & \multicolumn{1}{l|}{YES} & \multicolumn{1}{l|}{2.03e-9} & \multicolumn{1}{l|}{YES} & \multicolumn{1}{l}{6.00e-9} \\ \hline
\end{tabular}
\caption{Convergence behavior of the \texttt{chebop}, Fr\'echet, and $L$-scheme methods. The displayed residuals $\|\Psi(u)\|$\\ correspond to the last test case of the four infiltration problems.}
\label{tab:summary}
\end{table}

To conclude this section, we have seen that the easiest way to solve nonlinear one-dimensional BVP is provided by the \texttt{cheobop} class, when an initial guess of the solution which ensures the convergence is available. When the \texttt{cheobop} code fails to converge, the problem can be always solved by the explicit Newton-Fr\'echet linearization in function space and by the implicit $L$-scheme. These approaches generally require two or four times larger computing time, which however remains reasonable, of the order of tens of seconds. The alternative approaches are not only robust but also provide more computation details, as for instance graphs of the solution and update residuals. We also note that although it is slower than the Newton-Fr\'echet scheme, the $L$-scheme is much simpler, because it does not require the computation by hand of the Fr\'echet differentials.

\section{Exploring higher dimensions with Chebfun2 and Chebfun3}\label{sec:cheb23}

With the Cheffun2 \cite{TownsendandTrefethen2013} and Cheffun3 \cite{HashemiTrefethen2017} extensions, currently in the implementation stage, we may try to look for possible applications to infiltration problems in higher dimensions. 

Chebfub2 can solve linear partial differential equations defined on rectangular domains with the class \texttt{chebop2} \cite{TownsendandOlver2015}. However, a Newton-Fr\'echet method for unsteady infiltration in one spatial dimension cannot be implemented because only the determinant of the Jacobian but not the matrix itself can be constructed in Chebfun2 (see 'help jacobian' in Matlab command line). Trying to solve a linear problem with the implicit $L$-scheme by using \texttt{chebop2}, the code breaks down with the warning ``Warning: Matrix is close to singular or badly scaled'' and the error ``Solution was unresolved on a 65 by 65 grid''. The issue is not unexpected because typing ``help chebop2'' in the command line we get the information ``Warning: This PDE solver is an experimental new feature. It has not been publicly advertised''. As for Chebfun3, there are no solvers for linear problems.

The approach proposed in \cite{Driscolletal2008} uses Chebfun to solve steady-state problems in one spatial dimension and a FD scheme to perform the tine-stepping. In the following we refer to this approach as the FD-Chebfun coupling. Here, we propose an alternative Chebfun-FD coupling to solve unsteady problems in both one and two spatial dimensions, where the FD scheme solves steady-state problems at given time, while the time stepping is performed by Chebfun. A version of this coupling will be further used to evaluate orders of accuracy of the FD schemes.

Unlike in the steady-state problems solved in Section~\ref{sec:infil}, formulated in terms of pressure potential, now the unknown $\psi$ is the pressure head which is negative in the unsaturated regime. Following \cite{Suciuetal2021}, we consider the exact solution,
\begin{equation}\label{eq:exact_sol}
\psi^*(z,t) = - t\ z\ (z-1)\ -\ 0.5,
\end{equation}
and the constitutive laws for the water content $\theta$ and the conductivity $K$ given by
\begin{equation}\label{eq:param}
\theta(\psi) = \frac{1}{1-\psi}, \quad K(\theta(\psi)) = \psi^{2}.
\end{equation}

\subsection{FD-Chebfun coupling}\label{sec:ChebFDcheb}

The problem to be solved by Chebfun for the homogeneous, steady-state version of Richards equation (\ref{eq:Richards}) with initial and boundary conditions imposed by the exact solution (\ref{eq:exact_sol}) and parameterization (\ref{eq:param}) reads
\begin{align}
& \Phi(u)=u^{2}u''+2uu'(u'-1)=f,\label{eq:bvpR1}\\
& \alpha(u)\vert_{x=0}=u-u^*(0)=0,\label{eq:bvpR2}\\
& \beta(u)\vert_{x=1}=u-u^*(1)=0,\label{eq:bvpR3}
\end{align}
where $f=\Phi(u^*)$, $u=\psi(z,t)$, and $u^*=\psi^*(z,t)$. The FD-Chebfun coupling is achieved in the following by approximating the time derivative in (\ref{eq:Richards}) with the backward Euler discretization,

\begin{equation}\label{backwardEuler}
\frac{\partial}{\partial t}\theta(u)=\frac{1}{(1-u)^2}\frac{\partial u}{\partial t}
\approx\frac{1}{(1-u)^2}\frac{u(t)-u(t-\tau)}{\tau}.
\end{equation}
Equating the approximation (\ref{backwardEuler}) to $\Phi(u)$ (\ref{eq:bvpR1}) one obtains a new non-linear operator, 
\begin{equation}\label{eq:bvpR1t}
\widetilde{\Phi}(u)=-u+u_0+\tau(1-u)^2\left[\Phi(u)-f\right],
\end{equation}
where $u=u(t)$ and $u_0=u(t-\tau)$. The operator (\ref{eq:bvpR1t}) together with the boundary conditions (\ref{eq:bvpR2}) and (\ref{eq:bvpR3}) solve the unsteady problem with exact solution (\ref{eq:exact_sol}). The time stepping is initiated with $u_0=u^*(0)$ and is updated by $u_0=u$ at every step. 

\begin{figure}
\begin{minipage}[h]{0.45\linewidth}\centering
\includegraphics[width=\linewidth]{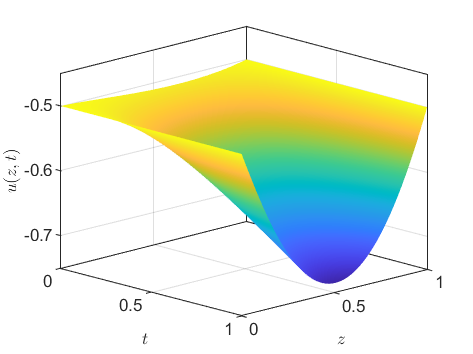}
\caption{\label{fig:chebR1D_sol.eps}Solution of FD-Chebfun coupling for unsteady Richards 1D equation.}
\end{minipage}
\hspace*{0.1in}
\begin{minipage}[h]{0.45\linewidth}\centering
\includegraphics[width=\linewidth]{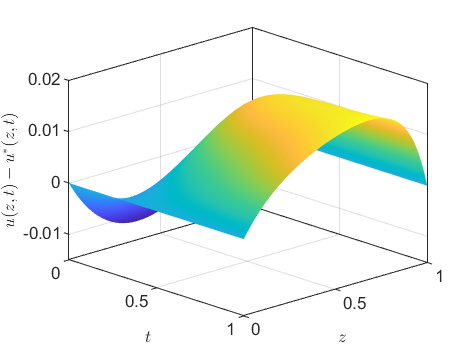}
\caption{\label{fig:chebR1D_err}Errors of FD-Chebfun coupling for unsteady Richards 1D equation.}
\end{minipage}
\end{figure}

The numerical and the exact solutions converted into \texttt{chebfun2} objects are used to compute errors (see code [\href{https://github.com/PMFlow/programming_with_chebfun/blob/main/Chebfun_FD/FD_Chebfun.m}{Chebfun\_FD/FD\_Chebfun.m}]). Figure~\ref{fig:chebR1D_sol.eps} represents the \texttt{chebfun2} solution after ten time steps, with $\tau=0.1$. The errors of order $10^{-2}$ presented in Fig.~\ref{fig:chebR1D_err} show that the overall accuracy of the coupled FD-Chebfun scheme is rather that of the FD scheme, not that of the Chebfun approximation (see also \cite{Driscolletal2008}). The code also computes residuals of the unsteady Richards equation which were found to be of the same order of magnitude as the errors shown in Fig.~\ref{fig:chebR1D_err}. 

\begin{figure}
\begin{minipage}[h]{0.45\linewidth}\centering
\includegraphics[width=\linewidth]{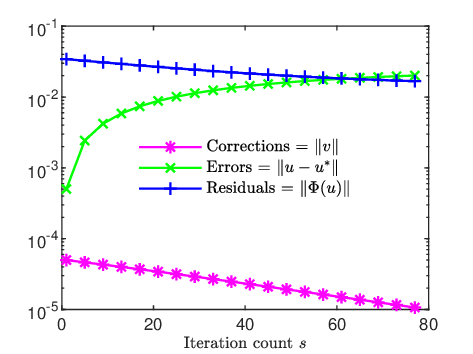}
\caption{\label{fig:conv_chebfun2FD}Convergence of the chebfun2-FD solution.}
\end{minipage}
\hspace*{0.1in}
\begin{minipage}[h]{0.45\linewidth}\centering
\includegraphics[width=\linewidth]{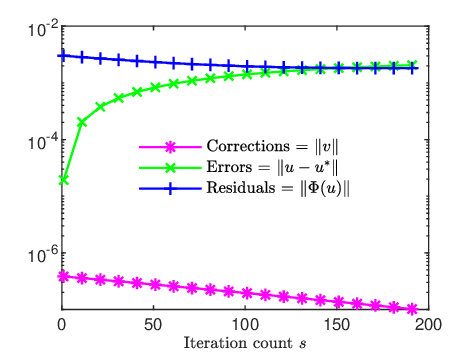}
\caption{\label{fig:conv_chebfun3FD}Convergence of the chebfun3-FD solution.}
\end{minipage}
\end{figure}

\subsection{Chebfun-FD coupling}\label{sec:ChebFD}

The coupling between Chebfun2 and Chebfun3 with one- and two-dimensional FD schemes is constructed as follows. At every time level, a steady state problem is solved by an explicit FD $L$-scheme. For the one-dimensional Richards equation (\ref{eq:Richards}) an explicit $L$-scheme is obtained by evaluating the second order operator in (\ref{eq:modPicard}) from the last iterate $\psi^{j,k}$,
\begin{equation}\label{eq:explL-Richards}
L(\psi^{j,k+1}-\psi^{j,k})-\frac{\partial}{\partial z}\left[K^{j,k}\frac{\partial}{\partial z}(\psi^{j,k}-z)\right]+\frac{\theta^{j,k}-\theta^{j-1}}{\tau}=0.
\end{equation}
The generalization of the scheme (\ref{eq:explL-Richards}) to higher dimensions is straightforward. Since the explicit $L$-schemes do not solve systems of algebraic equations, they can be up to ten times faster in solving unsteady problems (see \cite{Suciuetal2021} for implementation details and comparisons with an implicit scheme). In the following, the array representing the FD solution is imported into Chebfun2 (or Chebfun3) by interpolation procedures \cite{Trefethen2013,Trefethen2015}. The Chebfun code performs time and iteration steps, evaluates corrections, errors, and residuals and provide a new initial guess for the FD scheme. The coupling approach is summarized in Algorithm~\ref{alg:ChebFD}.

\begin{algorithm}\label{alg:ChebFD}
Two-dimensional Chebfun-FD coupling
\begin{verbatim}
# Define the domain and the independent variables:
D=[0 T y1 y2];
x = chebfun2(@(x,y) x,D); (time variable)
y = chebfun2(@(x,y) y,D); (space variable)
# Construct the exact solution:
E=chebfun2(@(x,y) -x.*y.*(1-y)-0.5,D);
# Differential operator and source term for Eq. (15):
phi=@(u) Op(u); f=phi(E);
# Define grid and initial guess for the FD scheme:
yy=y1:dy:y2; u0=E(0,:); p0=u0(yy); 
# Initialize variables:
t=0; k=1; correction(1)=1;
while t<=T
	while norm(correction(k)) > Tolerance
			[p,dt] = FD function (p0,f);
			# Convert FD solution to chebfun, compute corrections and errors:
			u=chebfun_interpolator(p,D);
			correction(k)=u-u0; u0=u;  
			error(k)=u-E(t,:); 
			# Compute unsteady FD solution, convert to chebfun2, compute residuals:
			q(:,time count)=p;
			u2=chebfun2_interpolator(q,D);
			residual(k)=phi(u2)-f;
			k=k+1; p0=p;
	end
t=t+dt; 
end
\end{verbatim}
\end{algorithm}

The code [\href{https://github.com/PMFlow/programming_with_chebfun/blob/main/Chebfun_FD/chebfun2FD.m}{Chebfun\_FD/chebfun2FD.m}] solves, according to Algorithm~\ref{alg:ChebFD}, the unsteady Richards equation (\ref{eq:Richards}) with initial and boundary conditions imposed by the exact solution (\ref{eq:exact_sol}) and the constitutive relations (\ref{eq:param}) in the domain $D=[0\; 0.1\;\; 0\; 1]$ in about 50 seconds. Figure~\ref{fig:conv_chebfun2FD} shows that while the corrections fall to the tolerance $\varepsilon=10^{-5}$, the errors and the residuals remain close to $10^{-2}$. Thus, the accuracy of the solution provided by the Chebfun-FD coupling is again much lower than the typical Chebfun accuracy.

In the case of two-dimensional unsteady version of (\ref{eq:Richards}), the exact solution (\ref{eq:exact_sol}) is replaced by $\psi^*(x,z,t) = - t\ x\ (x-1)\ z\ (z-1)\ -\ 0.5$ and the Chebfun3-FD solution is implemented in the code [\href{https://github.com/PMFlow/programming_with_chebfun/blob/main/Chebfun_FD/chebfun3FD.m}{Chebfun\_FD/chebfun3FD.m}]. Solving the problem in the domain $D=[0\; 0.03\;\; 0\; 1\;\; 0\; 1]$ takes a much larger time of about 25 minutes. Figure~\ref{fig:conv_chebfun3FD} shows that the errors and the residuals are again much larger than the corrections.

The same two-dimensional unsteady problem is solved with the FD L-scheme code in about 70 seconds [\href{https://github.com/PMFlow/programming_with_chebfun/blob/main/Chebfun_FD/FD2dimLscheme.m}{Chebfun\_FD/FD2dimLscheme.m}]. So, the coupling with Chebfun slows down the computation. Instead, importing the FD solution into Chebfun allows residuals evaluation. As seen in Figs~\ref{fig:conv_chebfun2FD}~and~\ref{fig:conv_chebfun3FD} the residuals converge to the error of the numerical solution with respect to a known exact solution. This suggests that the coupling Chebfun-FD can be used to evaluate the accuracy of the FD scheme, even when no exact solution is available, because the residuals can be always calculated by evaluating the differential operator on the numerical FD solution imported into Chebfun.

\subsection{Order of accuracy estimates}\label{sec:accuracy}

The accuracy of the classical discretization schemes is characterized by the so called observed order of accuracy \cite{Roy2005}. The order of accuracy $\alpha$ can be estimated by successively halving the mesh dimension $h$ according to
\begin{equation}\label{eq:accuracy}
\alpha=\ln\left(\frac{e_k}{e_{k+1}}\right)\;/\;\ln(2),\;\; k=1,2,\cdots,
\end{equation}
where $e_k=\|\psi_k-\psi^*\|$, is the error measured in the discrete  $L^2$ norm $\| \cdot \| := \sqrt{h} \| \cdot \|_{\ell^2}$, and $\| \cdot \|_{\ell^2}$ is the Euclidean vector norm. The procedure is described in Algorithm~\ref{alg:Accuracy}.

\begin{algorithm}\label{alg:Accuracy}
Accuracy estimates by two-dimensional Chebfun-FD coupling
\begin{verbatim}
# Define the domain and the independent variables:
D=[0 T y1 y2];
x = chebfun2(@(x,y) x,D); (time variable)
y = chebfun2(@(x,y) y,D); (space variable)
# Construct the exact solution:
E=chebfun2(@(x,y) -x.*y.*(1-y)-0.5,D);
# Differential operator and source term for Eq. (15):
phi=@(u) Op(u); f=phi(E);
for k=1: degree of the refinement
		[p,pE,q,qE] = FD function (k,f,E);
		# Compute the error of the FD solution 'p' w.r.t. the exact solution 'pE' at t=T: 		
		L2_p(k) = norm(p-pE);
		# Convert the time discrete dependent FD solution 'q' to chebfun2:
		u=chebfun2_interpolator(q,D);
		# Convert the discrete unsteady exact solution 'qE' to chebfun2:
		uE=chebfun2_interpolator(qE,D);
		# Compute the residual norm:
		R_p(k)=norm(phi(u)-phi(uE));
		# Use 'L2_p(k)' and 'R_p(k)' to compute orders of accuracy by Eq. (36): 
		- from errors of the discrete FD solution at the final time;
		- from residuals of the continuous space-time chebfun2 solution.
end
\end{verbatim}
\end{algorithm}

For the FD $L$-scheme, the mesh dimension is given by the constant space-step , $h=\Delta z$. The order of accuracy for the one-dimensional unsteady Richards equation estimated with the code [\href{https://github.com/PMFlow/programming_with_chebfun/blob/main/Accuracy/Accuracy_Chebfun2.m}{Accuracy/Accuracy\_Chebfun2.m}] according to (\ref{eq:accuracy}) from three successive halving of $\Delta z$ is shown in Table~\ref{table:acc_cheb2}. We remark that the estimates by errors of the FD $L$-scheme and by residuals of the \texttt{chebfun2} solution are very close to each other.
\begin{table}[h]
\begin{tabular}{ c c | c c c c c  c c c c}
  \hline
   & $\Delta z$ & & $\|\psi_k-\psi_{m}^*\|$ &  & $\alpha$ &  & $\|\Phi(u_k)\|$ &  & $\alpha$ &  \\
  \hline
   &  2.00e-01  & & 1.50e-02 &  & --   &  & 8.78e-02 &  & --   & \\
   &  1.00e-01  & & 3.61e-03 &  & 2.05 &  & 1.93e-02 &  & 2.18 & \\
   &  5.00e-02  & & 8.61e-04 &  & 2.07 &  & 4.57e-03 &  & 2.08 & \\
   &  2.50e-02  & & 1.32e-04 &  & 2.70 &  & 9.73e-04 &  & 2.23 & \\
  \hline
\end{tabular}
\caption{\label{table:acc_cheb2}Observed accuracy order $\alpha$ (\ref{eq:accuracy}) from errors at $t=T$ of the discrete FD solution\\ and residuals of the continuous space-time solution represented by a \texttt{chebfun2} object.}
\end{table}

Estimates of $\alpha$ for the two-dimensional unsteady Richards equation, computed from the same mesh refinement with the code [\href{https://github.com/PMFlow/programming_with_chebfun/blob/main/Accuracy/Accuracy_Chebfun3.m}{Accuracy/Accuracy\_Chebfun3.m}], are shown in Table~\ref{table:acc_cheb3}. The estimates by errors of the FD $L$-scheme and by residuals of the \texttt{chebfun3} solution are again very close to each other.

Since the computation of the residuals only uses the differential operator of the problem, the Chebfun-FD coupling provides a new method to evaluate the accuracy of the FD (and other classical discretization methods) without requiring the knowledge of the exact solution.

\begin{table}[h]
\begin{tabular}{ c c | c c c c c  c c c c}
  \hline
   & $\Delta z$ & & $\|\psi_k-\psi_{m}^*\|$ &  & $\alpha$ &  & $\|\Phi(u_k)\|$ &  & $\alpha$ &  \\
  \hline
   &  2.00e-01  & & 3.30e-02 &  & --   &  & 1.17e-01 &  & --   & \\
   &  1.00e-01  & & 7.27e-03 &  & 2.18 &  & 2.61e-02 &  & 2.16 & \\
   &  5.00e-02  & & 2.13e-03 &  & 1.77 &  & 7.70e-03 &  & 1.76 & \\
   &  2.50e-02  & & 1.35e-03 &  & 0.66 &  & 4.72e-03 &  & 0.71 & \\
  \hline
\end{tabular}
\caption{\label{table:acc_cheb3}Observed accuracy order $\alpha$ (\ref{eq:accuracy}) from errors at $t=T$ of the discrete FD solution\\ and residuals of the continuous space-time solution represented by a \texttt{chebfun3} object.}
\end{table}

\subsection{Failure of modeling saturated-unsaturated transition in Chebfun}\label{failure}

\begin{figure}
\begin{minipage}[h]{0.45\linewidth}\centering
\includegraphics[width=\linewidth]{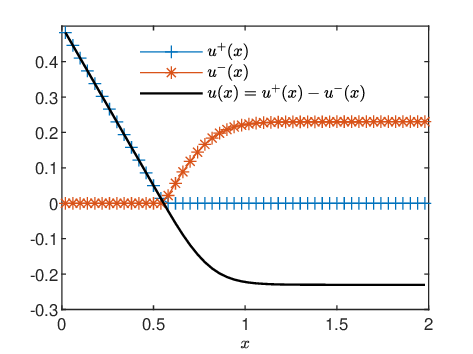}
\caption{\label{fig:u_pm}Splitting into negative and positive parts $u=u^-+u^+$.}
\end{minipage}
\hspace*{0.1in}
\begin{minipage}[h]{0.45\linewidth}\centering
\includegraphics[width=\linewidth]{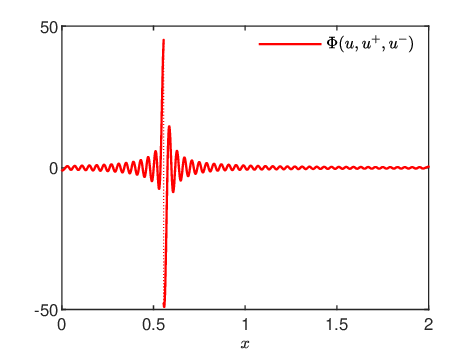}
\caption{\label{fig:Phi_pm}Solution of the splitting operator $\Phi(u,u^+,u^-)$.}
\end{minipage}
\end{figure}

The Chebfun-FD coupling is illustrated in Sections~\ref{sec:ChebFD}~and \ref{sec:accuracy} for unsaturated flow regime. The natural question is, does the coupling also work for saturated-unsaturated transitions? Since in saturated regime the conductivity $K$ is constant and Richards equation degenerates into an elliptic equation, we can try the splitting
$\Phi=\Phi^-+\Phi^+$, where $\Phi^-$ and $\Phi^+$ are differential operators corresponding to the parabolic Richards equation, operating in the unsaturated regime, and to the elliptic equation in saturated regime. Chebfun operators operate on the whole function $u$, unlike the finite difference operators which automatically change the form for 
$\psi\geq 0$ and $\psi<0$. Therefore, we have to split $\psi$ into its negative and positive part, $\psi^-=\max(0,-\psi)$, $\psi^+=\max(0,\psi)$, and export them as chebfun objects, $u^-$ and $u^+$, as shown in Fig.~\ref{fig:u_pm}. With these, the splitting operator becomes
\begin{equation}\label{eq:splitting}
\Phi(u,u^+,u^-)=\frac{1-\sign(u)}{2}\Phi^-(u^-)+\frac{1+\sign(u)}{2}\Phi^+(u^+),
\end{equation}
where the $\sign$ function prevents spurious oscillations that may appear when $\Phi^-$ operates for $u\geq 0$ and  $\Phi^+$ operates for $u<0$ (see Fig.~\ref{fig:u_pm}).

The problem for the one-dimensional unsteady Richards equation is implemented in [\href{https://github.com/PMFlow/programming_with_chebfun/blob/main/Accuracy/Failure.m}{Accuracy/Failure.m}] by the following code lines:
\begin{verbatim}
[p] = FDL_scheme(E,I);        % FD solution for exact solution E and grid dimension I.
u=chebfun(p',D,'equi');
pm=max(0,-p);                 % negative part of p (unsaturated)
pp=max(0,p);                  % positive part of p (saturated)
um=chebfun(pm',D,'equi');
up=chebfun(pp',D,'equi');
phi_m = @(u) diff(u,2)+c*diff(u).*(diff(u)+1); % operator for unsaturated regime
phi_p = @(u) diff(u,2);                        % operator for saturated regime
phi = @(u,up,um) (1-sign(u))/2.*phi_m(-um)+(1+sign(u))/2.*phi_p(up);
residual=norm(phi(u,up,um));
\end{verbatim}

\begin{figure}
\begin{minipage}[h]{0.45\linewidth}\centering
\includegraphics[width=\linewidth]{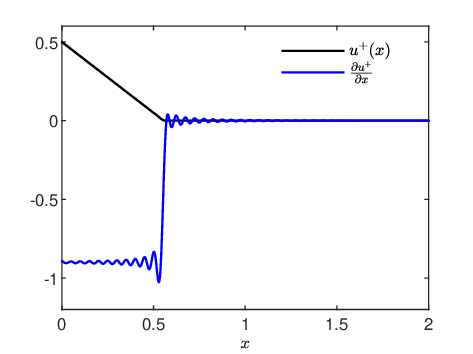}
\caption{\label{fig:du_p}First order derivative of $u^+$.}
\end{minipage}
\hspace*{0.1in}
\begin{minipage}[h]{0.45\linewidth}\centering
\includegraphics[width=\linewidth]{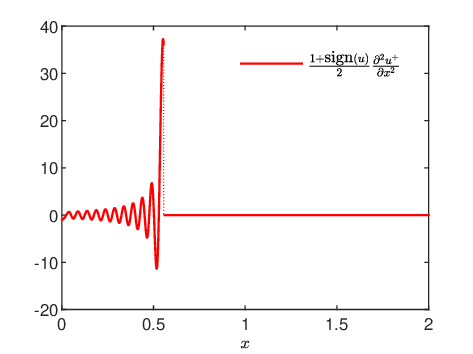}
\caption{\label{fig:d2u_p}Second order derivative of $u^+$.}
\end{minipage}
\end{figure}
The code operates on the steady-state solution \cite[Scenario 1]{Suciuetal2021} which solves the Richards equation with the Gardner parameterization. While the FD $L$-scheme verifies the solution with errors of the order $10^{-7}$, the splitting fails and produces the highly oscillatory solution shown in Fig.~\ref{fig:Phi_pm}. The issue can be investigated by analyzing the effect of the differentiation on the positive part of the solution, $u^+$. The first derivative shown in Fig.~\ref{fig:du_p} already oscillates close to the inflection point of the solution $u$. Figure~\ref{fig:d2u_p} shows a dramatic increase of the oscillations of the second order derivative, which coincides with the left side of the solution of the splitting operator shown in Fig.~\ref{fig:Phi_pm}, corresponding to the saturated regime. Hence, we found that Chebfun successfully resolves the solution $u$ but fails to resolve its derivatives.

Summarizing this section, we have seen that the Chebfun-FD coupling is a valuable tool for accuracy assessments of the classical discretization schemes for Richards equation, but it only works for the unsaturated flow regime.

\section{conclusions}\label{sec:concl}
The main results of this study are summarized as follows.

$\bullet$ The automated \texttt{chbeop} class of Chebfun solves nonlinear steady-state BVP for Richards equation in one spatial dimension to roughly machine precision, provided that a suitable initial guess which ensures the convergence of the Newton iterations is available.

$\bullet$ If \texttt{chbeop} fails, a well posed BVP can always be solved, within a wider range of initial iterates, by an explicit linearization of the nonlinear operator in function space based on Fr\'echet differential and the linear Chebfun solver, as well as by globally convergent $L$-schemes, quasi-Newton methods where the derivatives are replaced by suitably chosen positive numbers.

$\bullet$ The coupling of Chebfun2 and Chebfun3 with explicit $L$-schemes can be used to solve unsteady infiltration problems in one and two spatial dimensions and is mainly useful in assessing the accuracy of the classical discretization schemes.

$\bullet$ Because Chebfun fails to resolve the derivatives of the solution close to the inflection point, where the pressure head changes the sign, Chebfun-based applications are only possible for the unsaturated flow regime.

\end{document}